\documentclass{article}
\usepackage[ansinew]{inputenc}
\usepackage[english]{babel}
\usepackage{latexsym}
\usepackage{amsthm}
\usepackage{amsmath}
\usepackage{amsfonts}
\usepackage{amssymb}
\usepackage{graphicx}
\usepackage{times}
\usepackage{url}
\usepackage{float}

\begin{document}
\newcommand{\suces}[3]{{#1_{#2},\hdots,#1_{#3}}}
\newcommand{\realamp}{\R\cup\{-\infty\}}
\newcommand{\A}{\mathbb{A}}
\renewcommand{\P}{\mathbb{P}}
\newcommand{\T}{\mathbb{T}}
\newcommand{\ptoo}[1]{{(#1_1,#1_2)}}
\newcommand{\pto}[2]{{(#1,#2)}}
\newcommand{\ptoop}[1]{{[#1_1,#1_2,#1_3]}}
\newcommand{\ptop}[3]{[#1,#2,#3]}

\newcommand{\m}{\medskip}

\newcommand{\N}{\mathbb{N}}
\newcommand{\Q}{\mathbb{Q}}
\newcommand{\R}{\mathbb{R}}

\newcommand{\CC}{\mathcal{C}}
\newcommand{\HH}{\mathcal{H}}
\newcommand{\TT}{\mathcal{T}}

\newcommand{\id}{\operatorname{id}}
\newcommand{\im}{\operatorname{im}}

\newcommand{\V}{\operatorname{\tt{V}}}

\newcommand{\est}{\operatorname{est}}
\newcommand{\conv}{\operatorname{conv}}

\newtheorem{thm}{Theorem}
\newtheorem{lem}{Lemma}
\newtheorem{dfn}{Definition}
\newtheorem{rem}{Remark}
\newtheorem{ex}{Example}
\newtheorem{cor}{Corollary}

\title{A note on  tropical triangles in the plane}
\author{M. Ansola and M.J. de la
Puente\thanks{Departamento de Algebra, Facultad de Matem{\'a}ticas,
Universidad Complutense, 28040--Madrid, Spain}\ \thanks{Partially
supported by MTM 2005--02865 and UCM 910444; e--mail: mpuente@mat.ucm.es.}} \maketitle

2000 Math. Subj. Class.: 52C35, 52C20, 15A39,   12K99.


\begin{abstract}
We define transversal tropical triangles  (affine and projective) and
characterize them via six inequalities to be satisfied by
the coordinates of the vertices. We prove that the vertices of a  transversal tropical triangle
are tropically independent and they tropically span a classical hexagon whose sides have slopes $\infty,0,1$.
Using this classical hexagon, we determine a parameter space for transversal tropical triangles.
The coordinates of the vertices of a transversal tropical triangle determine a
tropically regular matrix. Triangulations of the tropical plane are obtained.
\end{abstract}

\section{Introduction}
Triangles are, after points
and lines,  the simplest figures in any geometry.
They can be defined either by three different non--collinear points, called vertices,
or by three non concurrent pairwise transversal lines, called sides. The vertices of a triangle span the plane
(affine or projective) and  they are independent points.
Also, triangles are the simplest two--dimensional convex figures and they provide  tillings  of the plane, or triangulations.
This is all elementary mathematics.

\m
 If we move to the tropical plane, then  natural questions  about triangles arise.  How are the notions of span, independence,  convexity and transversality (in their tropical versions) related to tropical triangles? On the one hand, in \cite{Cuninghame_B}, \cite{Wagneur_M}  it is shown that no finite family of points can tropically span the plane.
On the other hand, tropical convexity has been thoroughly studied in \cite{Develin}; in particular, tropical triangles are defined there,
and five combinatorial types of tropical triangles are shown to exist, up to symmetry.
However, most  such tropical triangles have non--transversal  sides.

\m
In this paper, we propose a finer definition of tropical triangle. 
Such tropical triangles will be called \emph{transversal}. It amounts to restricting to  one combinatorial type from \cite{Develin}.  Our starting point (from elementary geometry) is the following: the  tropical triangles we are interested in  are only those defined by \emph{three different non--collinear points  $a,b,c$ such
that, when joined by pairs yield three different lines $ab,bc,ca$ which, when
intersected by pairs return the original points}. 
The  intersection (resp. join) we are talking about here  is \emph{stable intersection} (resp. \emph{stable join}) and  these are the right notions to consider in tropical geometry. And it turns out that, in transversal conditions  stable intersection (resp. stable join) is nothing but plain intersection  (resp. join).

\m
It is a basic fact  that a tropical line in the plane carry a special point, called \emph{vertex}. Therefore, to a tropical triangle $T$ we can associate a family of six points: the three vertices of $T$ and the three vertices of the tropical sides of $T$.
In theorem \ref{thm:hexagono} it is proved that the tropical span of the vertices  of a transversal tropical
triangle $T$ equals the classical convex hull of the six related points. This means that giving a  transversal tropical triangle $T$ amounts to giving a classical hexagon, $\HH(T)$, the sides of which have slopes  $\infty,0,1$.
It is our opinion that all triangles should look alike, in any geometry and this is
not the case for  tropical triangles, as defined 
in \cite{Develin},
but it is certainly true, for
transversal tropical triangles. Moreover, the lattice lengths the sides of $\HH(T)$  parameterize the tropical triangle $T$.

\begin{figure}[ht]
 \centering
  \includegraphics[width=10cm]{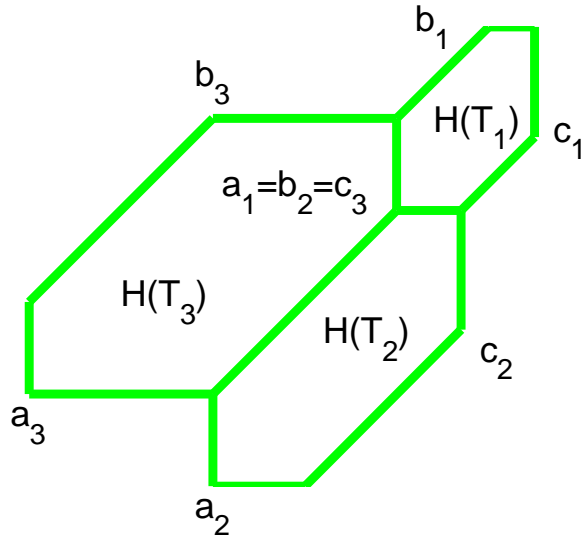}\\
  \caption{Tropical triangulation of the plane.}
  \label{fig:trian_1}
  \end{figure}

\m
In \cite{Joswig}, Joswig rises the question of what should be
the right notion of \emph{tropical triangulation}. For the
tropical plane (affine or projective), we give the
following solution: a triangulation $\TT$ of the tropical plane
 is a family of transversal tropical triangles  $\{T_j:j\in J\}$ such
that the  associated family
$\{\HH(T_j):j\in J\}$ tessellate the classical  plane. More precisely, if
two tropical triangles $T_1,T_2$ in $\TT$ meet, all
they share is one vertex and one side and, moreover, the
associated   classical hexagons $\HH(T_1),\HH(T_2)$ have just one
side in common, including the two end points, see figure \ref{fig:trian_1}.

\m
The main results of the paper are theorems \ref{thm:triangle_affine} and
\ref{thm:triangle_projective}, where
transversal tropical triangles  are  characterized by six strict inequalities to be satisfied by the coordinates of the vertices (cases affine and projective).  The inequalities in theorem \ref{thm:triangle_projective} show  a high level of symmetry.

\m
We  work exclusively in the tropical setting. Some authors solve problems in
tropical geometry by the \emph{lifting method}, see \cite{Izhakian}, \cite{Richter}, \cite{Tabera_Pap}.
This means that they start with a question in tropical
geometry, lift the question to classical  geometry (when
possible, see \cite{Jensen}, \cite{Tabera_Trans}), solve the question there (when possible) and
tropicalize the solution (always possible). Although the
existence of liftings is one powerful reason to do
tropical geometry, we will not use liftings at all.

\m
\emph{Tropical geometry} (and the study of certain closely related objects called \emph{amoebas}) is a very recent trend in Mathematics, see
\cite{Ansola}, \cite{Einsiedler}, \cite{Gathmann}, \cite{Gelfand}, \cite{Itenberg}, \cite{Mikhalkin_W}, \cite{Richter}, \cite{Sturmfels}, \cite{Viro_D}, \cite{Viro_W}. Its algebraic counterpart is \emph{tropical algebra} and it has been studies since about 1950. It is
related to control theory, automata theory, scheduling
theory, discrete event systems, optimization, combinatorics, mathematical
physics, etc. It has applications in
complex and real enumerative geometry, phylogenetics, etc., see \cite{Gathmann_M}, \cite{Mikhalkin_E}, \cite{Mikhalkin_T}, \cite{Shustin}. It is a fast developing geometry, following the track of algebraic geometry, see \cite{Mikhalkin_J}, \cite{Speyer}, \cite{Vigeland}. It is connected with toric geometry, see \cite{Hsie}, \cite{Nishinou}, \cite{Shustin_W}.

\m
Tropical algebra has appeared in the literature
under various denominations such as minimax--algebra,
max--algebra, min--algebra,  max--plus algebra, min--plus algebra,  semirings, modulo\"{\i}ds, dio\"{\i}ds, pseudorings,
pseudomodules, band spaces over belts, idempotent mathematics (semirings,
analysis, calculus, etc.), Maslov dequantization, etc., see
\cite{Akian}, \cite{Baccelli}, \cite{Butkovic}, \cite{Cuninghame}, \cite{Gaubert}, \cite{Litvinov}, \cite{Viro_D}, \cite{Wagneur_F}.  

\m
As a rule, we will  use the adjective classical (for classical mathematics) as opposed to tropical.
In this note we present results which can be traced back to  \cite{Ansola},
 but have been very much elaborated afterwards. Transversal triangles are called stable triangles there.
 A previous version of this paper can be found in ArXiv.

\m
We take the opportunity to  thank B. Bertrand for his
excellent introductory talks about  Tropical Geometry,
held  in the academic course 2004--05 at Departamento de
Algebra, Facultad de Matem\'{a}ticas, Universidad
Complutense de Madrid (Spain) and to L.F. Tabera for his
interest.

%
%
\section{Notations and background on elementary tropical geometry}

The \emph{tropical semi--field} is the set $\T:=\realamp$ endowed with
tropical addition $\oplus$ and tropical multiplication $\odot$.
These operations are defined as follows:
    $$a\oplus b=\max\{a,b\},\qquad a\odot b=a+b,$$ for
    $a,b\in\realamp$. 
    Note that tropical addition is \emph{idempotent},
    i.e.,
    $a\oplus a=a$, for $a\in\T$. Tropical addition  is
    associative,
    commutative and  $-\infty$ is the neutral element.
    Tropical multiplication  is associative,
    commutative and  $0$ is the neutral element. The
    element $-a$ is  inverse to $a$ with respect to $\odot$,
     for    $a\in \R$. Moreover, multiplication is distributive
     over addition, since
        $$a+\max\{b,c\}=\max\{ a+b, a+c \}.$$
However, we
    cannot find an inverse element,
    with respect to $\oplus$, for any $a\in\R$ and
    this is why $\T$ is NOT  a field.

\m
For $n\in \N$, the  tropical affine $n$--space is
$\T^n$, where tropical addition and multiplication are defined
coordinatewise.  For addition, the neutral element is $p_n:=(-\infty, 
\ldots,-\infty)$, but one must realize that     $a\oplus b=b$ does not imply  $a=p_n$.

\m
The
tropical projective $n$--space, $\T\P^n$, is the defined
as follows. In the space
$\T^{n+1}\setminus\{p_{n+1}\}$ we define an
equivalence relation $\sim $ by letting
$(b_1,\ldots,b_{n+1})\sim (c_1,\ldots,c_{n+1})$ if there
exist $\lambda\in\R$ such that
    $(b_1+\lambda,\ldots,b_{n+1}+\lambda)=
    (c_1,\ldots,c_{n+1}).$
The equivalence class of $(b_1,\ldots,b_{n+1})$ is denoted
$[b_1,\ldots,b_{n+1}]$ and its elements are obtained by
adding  multiples of the vector $(1,\ldots,1)$ to the point
$(b_1,\ldots,b_{n+1})$.

\m Points in $\T^n$ (resp. $\T\P^n$) without infinite
coordinates will be called \emph{interior points}. The
rest of the points will be called \emph{boundary points}.
The \emph{boundary} of $\T^n$ (resp. $\T\P^n$)
 is the union of its boundary points. Note that
 $p_n$ is a boundary point in $\T^n$. 

\m In this paper we only work  in the tropical plane
(affine or projective). Therefore,  we will just present
the following notions for $n=2$, although  they  apply in
any dimension $n\ge2$. We will use $X,Y,Z$ as variables.

The tropical projective plane $\T\P^2$ is covered by three copies
of the tropical affine plane $\T^2$ as follows. The maps
$$j_3:\T^2\to\T\P^2,\quad (x,y)\mapsto \ptop xy0,\qquad j_2:\T^2\to\T\P^2,\quad (x,z)\mapsto \ptop x0z,$$
    $$j_1:\T^2\to\T\P^2,\qquad (y,z)\mapsto \ptop 0yz$$
    are injective and we have
$\T\P^2=\im j_3\cup\im j_2\cup\im j_1.$ The
complementary set  of $\im j_3$ is
$$\T\P^2\setminus\im j_3=\{\ptop
xy{-\infty}:x,y\in\T\}$$ and it is in bijection with
$\T\P^1$ (forget the last coordinate!) 
Similarly for $j_k$, $k=1,2$. It is easy to check that the
set of interior points
 of $\T\P^2$ is equal to the intersection
$\im j_3\cap \im j_2\cap \im j_1.$


\m As we already know, the projective tropical coordinates
of a point in $\T\P^2$ are not unique. In order to work
with a unique triple of coordinates for (almost) each
point, we choose a \emph{normalization}. A few points in
$\T\P^2$ do not admit normalized coordinates, but this
will not be a serious obstacle. Our favorite normalization
is  making the last coordinate equal to zero. We  will express this by  saying that
\emph{we work in} $Z=0$; it means passing from the
projective plane to the affine one, via $j_3$. Of course, the point
$\ptop ab{-\infty}$ does not admit normalized coordinates
in $Z=0$, for  $a,b\in\T$. Other possible normalizations
are $Y=0$, or $X=0$, or $X+Y+Z=0$, etc.

\m
Let $n\in \N$. We can consider
\emph{tropical polynomials} in any number of variables
$\suces X1n$ with coefficients in $\T$. Write $\overline{X}$ for $(\suces X1n)$ and
let $i=(\suces i1n)\in \N^n$ be a multi--index. Then write
    $$\overline{X}^{\odot i}=X_1^{\odot {i_1}}\odot\cdots\odot
      X_n^{\odot {i_n}}=i_1X_1+\cdots  + i_nX_n$$ and let
        $$p(\overline{X})=\bigoplus_{i\in I}a_i\odot \overline{X}^{\odot i}=
        \max_{i\in I}\{a_i + i_1X_1+\cdots  + i_nX_n\},$$
where $I\subset \N^n$ is some finite set and $a_i\in\T$. Being
$-\infty$  the neutral element for tropical addition,
 terms  having $a_i=-\infty$ may be omitted in $p$.
The polynomial $p$  is \emph{homogeneous} if there exists
$d\in\N$ such that $i_1+\cdots+ i_n=d$, for all $i\in I$ with
$a_i\neq-\infty$. If $p$ is homogeneous and
$a_i\neq-\infty$, FOR ALL
    $$i\in\{(d,0,\ldots,0),(0,d,0,\ldots,0),\ldots,
    (0,\ldots,0,d)\},$$ then we say that $p$  has
    \emph{degree $d$}. By means of an extra variable, we can
\emph{homogenize} a non--homogeneous tropical polynomial
$p$, easily. We will use capital letters to denote
homogeneous polynomials and small letters to  denote
arbitrary polynomials. The non--homogeneous polynomial $p$
is said to have \emph{degree $d$}
 if its homogenization has degree $d$.  In particular,
  the degree is NOT defined for some tropical
polynomials. 
In the literature one can find a more general notion of tropical degree, but this one is good enough for our purposes.

\m A tropical polynomial  $p$ (resp. homogeneous polynomial $P$) of
degree $d>0$ in $n$ (resp. $n+1$) variables defines a  so
called \emph{tropical affine hypersurface $\V(p)$ (resp. projective)
hypersurface $\V(P)$}
 in $\T^n$ (resp. $\T\P^n$). By definition, $\V(p)$ (resp.
$\V(P)$) is the set of \emph{points in $\T^n$ (resp.
$\T\P^n$) where the maximum is attained, at least, twice}.
This is certainly DIFFERENT from the classical  definition of
algebraic hypersurface. A tropical hyperplane is a tropical hypersurface defined by a linear polynomial.


\m Let  $n=2$. In this case, we use variables $X,Y,Z$,
instead of $X_1,X_2,X_3$. In the tropical plane,
hypersurfaces are called \emph{tropical curves}. We have
\emph{tropical lines, conics, cubics}, etc., meaning
curves defined by  tropical polynomials of degree 1, 2, 3,
etc., which are homogeneous in the projective case. Every tropical projective curve $\CC$ is covered by
three associated affine curves, namely,
$\CC=\CC_3\cup\CC_2\cup\CC_1,$ where $\CC_k:=\CC\cap \im j_k$, $k=1,2,3$.

\m The simplest tropical plane curves are lines, of course. A
\emph{tropical line} in the affine plane  is
$\V(p)$,
 where
    $$p=a\odot X\oplus b\odot Y\oplus c=\max\{a+X ,b+Y,c
     \}$$
for some coefficients $a,b,c\in\R$. Notice that the point
$p_2=(-\infty,-\infty)\in\T^2$ does not belong to $\V(p)$. 
The line $L:=\V(p)$  is easy to describe. The points  $\pto {c-a}{c-b}$, $\pto
{-\infty}{c-b}$ and  $\pto {c-a}{-\infty}$ belong to $L$; the first one is interior, while the other two are boundary points.
Moreover, $L$ is the union of three  rays,
meeting at the interior point $\pto {c-a}{c-b}$.  The directions of these rays
are West,   South and North--East. Notice that there is no boundary point at the end of the North--East ray of $L$.
The homogenization of $p$ is
    $$P=a\odot X\oplus b\odot Y\oplus c\odot Z=\max\{a+X ,
    b+Y,c+Z
     \}.$$
$P$ defines the \emph{tropical line} $\overline{L}:=\V(P)$  in $\T\P^2$ and $L$ embeds in $\overline{L}$ via
$j_3$. In particular, we have
$$j_3\pto{c-a}{c-b}=\ptop{c-a}{c-b}0=\ptop{-a}{-b}{-c},$$
and this point is called the \emph{vertex} of $\overline{L}$.
We also have
$$j_3\pto{c-a}{-\infty}=\ptop{-a}{-\infty}{-c},\qquad j_3\pto{-\infty}{c-b}=\ptop{-\infty}{-b}{-c},$$
and we find that $\ptop {-a}{-b}{-\infty}$ is the only point in
 $\overline{L}\setminus L$. Actually, this is the  missing boundary point in the North--East ray of $L$.
In addition, the identification of $L$  with its image in $\overline{L}$ via $j_3$ allows us to have a graphical
representation of $\overline{L}$ in $Z=0$, where the only missing point is  $\ptop {-a}{-b}{-\infty}$.
Of course, we can also represent
$\overline{L}$ in $Y=0$ or in $X=0$, easily.

\m
Let $n\in\N$. Given an $n\times n$ matrix $A$ with real
entries, the \emph{tropical determinant} of $A$, (also called \emph{permanent}) is defined as
follows:
$$|A|_{trop}:=\left|\begin{array}{ccc}
a_{11}&\cdots&a_{1n}\\
\vdots&&\vdots\\
a_{n1}&\cdots&a_{nn}\\
\end{array}\right|_{trop}=\max_{\sigma\in S_n}\{a_{1\sigma(1)}+\cdots+a_{n\sigma(n)}\}$$
where $S_n$ denotes the permutation group in $n$ symbols.
The matrix  $A$ is \emph{tropically singular} if the
maximum in $|A|_{trop}$ is attained, at least,
twice. Otherwise, $A$ is \emph{tropically regular}, (also said that $A$ has a \emph{strong permanent}).


\m
There exists a \emph{duality} between lines in the tropical
projective plane and interior points in the projective plane. Given an
interior point $a=\ptoop a$ in $\T\P^2$, let $L_a$ denote
the line in $\T\P^2$ 
defined by tropical linear form
 is $a_1\odot X \oplus a_2\odot Y
\oplus a_3\odot Z.$ Obviously, we have
$$b\in L_a\iff a\in L_b,$$  meaning that
$\max\{a_1+b_1,a_2+b_2,a_3+b_3\}$ is attained, at least, twice.

\m Let two points $a,b$ in the tropical plane (affine or projective) be given. If $a$ and $b$ do not both lie on a classical line of slope $\infty,0,1$, then there exists a unique tropical line through both points and this line is called the (tropical) \emph{join} of $a$ and $b$. Otherwise, there  exist infinitely many
tropical lines going through $a$ and $b$.  The (tropical)
\emph{stable join}   of $a,b$  is defined as the limit, as
$\epsilon$ tends to zero, of the tropical lines going
through perturbed points $a^{v_\epsilon}, b^{v_\epsilon}$.
Here, $a^{v_\epsilon}$ denotes a translation of $a$ by a
length--$\epsilon$ vector  $v_\epsilon$,  see
\cite{Gathmann},\cite{Richter}. We denote this
line by $ab$. Of course, if $a$ and $b$ do not both lie on a classical line of slope $\infty,0,1$,
then their join and their stable join coincide. 

\m Now, the intersection of  two tropical lines $L,M$ in the  plane (affine or projective)  may be a  point or a ray.
In the former case, we will say that the tropical \emph{lines $L$ and $M$ are transversal}.
In any case, we define the so called
\emph{stable intersection}, $L\cap_{\est}
M$, as the limit point, as $\epsilon$
tends to zero, of the intersection of perturbed lines
$L^{v_\epsilon}, M^{v_\epsilon}$. Here,
$L^{v_\epsilon}$ denotes a translation of $L$ by a
length--$\epsilon$  vector $v_\epsilon$.

\m It is well known that duality  transforms
stable join  into stable intersection and conversely,
 i.e.,  $$L_a\cap_{\est}L_b=c\iff ab=L_c,$$ for
 $a,b,c$ interior points in $\T\P^2$.  By duality, we will say that the points
 $a,b$ are \emph{transversal} if there exists a unique tropical line passing $a$ and $b$.

\m Stable intersection and stable join are defined in wider generality.
Now consider $n$ hyperplanes in $\T\P^n$ and
take an associated system of $n$ linear tropical homogeneous
polynomials in $n+1$ variables. Let $A$ be the $n\times
(n+1)$ coefficient matrix of the system. For each
$j=1,\ldots, n+1$, let $A^j$ be the square matrix obtained  by deleting the $j$--th column
from $A$. Then the
\emph{tropical version of Cramer's rule} tells us that the point
$[|A^1|_{trop}, \ldots,|A^{n+1}|_{trop}]\in\T\P^n$ is the stable
intersection of the $n$ given  hyperplanes. Moreover, the intersection of the $n$ hyperplanes
equals the stable intersection if and only if  $A^j$ is tropically regular, for all $j=1,\ldots, n+1$; see
\cite{Richter}.

 \m  For tropical lines in the plane, the tropical version of Cramer's rule  goes as follows:
 the stable
intersection of the lines $L_a$ and $L_b$ is the point
$$[\max\{a_2+b_3,b_2+a_3\},
 \max\{a_1+b_3,b_1+a_3\},
    \max\{a_1+b_2,b_1+a_2\}].$$ Since the computation of this point is nothing but
      a tropical version of the
    cross--product of the triples $a$ and $b$,  we will
     denote it by $a\otimes b$. Notice that $a\otimes b=b\otimes a$.
     In other words,  the tropical version of Cramer's rule  in the plane means
     $$L_a\cap_{est}L_b=a\otimes b,$$
and, by duality,
$$ ab=L_{a\otimes b}.$$
In particular, the vertex  of the tropical line $ab$ is the point $-(a\otimes b)$.

\m Tropical cross--product satisfies $a\otimes a=-a$, if $a$ is interior. But, unluckily, tropical cross--product is nonassociative and so,  it seems tricky to compute expressions such as $(c\otimes a)\otimes (a\otimes b)$. However, we will see  in corollary \ref{cor:tres_val} that only three values are possible for this long expression.

\begin{lem}\label{lem:transversal}
Interior points $a,b\in\T\P^2$ are transversal if and only if $a\otimes b\not\in\{-a,-b\}$.
\end{lem}
\begin{proof}
If $a,b$ are interior points in $\T^2$, then a simple computation shows that $a\otimes b\in\{-a,-b\}$ if and only if the points $a,b$ lie on a classical line of slope $\infty,0,1$, and this is the non--transversal case.
Now if $a,b$ are interior points in $\T\P^2$, we obtain the desired result, either passing to the affine setting, or using Cramer's rule. Indeed, an easy computation shows that $a\otimes b\in\{-a,-b\}$ if and only if
some coordinate in $a\otimes b=[\max\{a_2+b_3,b_2+a_3\},
 \max\{a_1+b_3,b_1+a_3\},
    \max\{a_1+b_2,b_1+a_2\}]$ is tropically singular.
\end{proof}

\begin{lem}\label{lem:long_expr}
Suppose that $a,b,c$ are three interior points in the tropical plane (affine or projective) which are tropically non--collinear. Assume that
$a,b$ are transversal, $a,c$ are transversal and the tropical lines $ab, ca$ are transversal. Then  $(c\otimes a)\otimes (a\otimes b)=a$, i.e.,  $ca\cap_{est}ab=ca\cap ab=a$.
\end{lem}

\begin{proof} Since the lines $ab$ and $ca$ are transversal,  then the points
   $c\otimes a,  a\otimes b$ do not lie on a classical line of slope $\infty,0,1$. Since the points $a,b$ are transversal, then  $a$ and $b$ lie on different rays of the tropical line with vertex at $-(a\otimes b)$. A similar situation is true for the points $a,c$ and the tropical line  with vertex at $-(c\otimes a)$. Then $a$ is the unique point of intersection of $ab$ and $ca$. 
This gives a picture in $\R^2$. A symmetric picture, with respect to the origin, is obtained by considering the points $c\otimes a$, $a\otimes b$ and the vertex, $-((c\otimes a)\otimes (a\otimes b))$, of the unique tropical line through them. Therefore  $(c\otimes a)\otimes (a\otimes b)=a$.
\end{proof}

The corollary below is a direct consequence of lemmas \ref{lem:transversal} and  \ref{lem:long_expr}.
\begin{cor}\label{cor:tres_val}
Suppose that $a,b,c$ are three interior points in the tropical plane (affine or projective) which are tropically non--collinear. Then
$$(c\otimes a)\otimes (a\otimes b)\in\{-(c\otimes a), -(a\otimes b), a\}.\qed$$
\end{cor}

\section{Transversal triangles in the tropical plane}
In this section, $a,b,c$ will always denote three different interior points in the tropical plane (affine or projective) which are tropically non--collinear. The  \emph{tropical sides} defined by $a,b,c$ are the
tropical lines $ab, bc, ca$ and
we know that the points $-(a\otimes b)$, $-(b\otimes c)$ and $-(c\otimes a)$ are the  vertices of them.

\m
In the first example below, we see  that
$ab\cap_{est}bc\neq b$ and this is unpleasant for a triangle.
In the second one,
we see  have $ca\cap_{est}ab= a, ab\cap_{est}bc= b, bc\cap_{est}ca= c$, showing the  kind of triangles we are interested in.


\begin{ex}\label{ex:non_good}
In $\T\P^2$  take $a=\ptop {-1}10$, $b=\ptop 000$,
$c=\ptop {-1}20$.  Then, $a\otimes b=[1, 0, 1]$ and $b\otimes c=[2, 0, 2]$, so that
$(a\otimes b)\otimes(b\otimes c)=[2,3,2]=-(a\otimes b)\neq b$.
\end{ex}

\begin{ex} \label{ex:good_and_prop}
The reader can easily check that the points $a=\ptop {-3}{-1}0$, $b=\ptop 000$, $c=\ptop {-1}20$ satisfy
$(c\otimes a)\otimes(a\otimes b)=a$, $(a\otimes b)\otimes(b\otimes c)= b$ and $(b\otimes c)\otimes(c\otimes a)= c$.
\end{ex}

Thus,  for some purposes, some care must be taken in order to define triangles
in the tropical plane (affine or projective).

\begin{figure}[ht]
 \centering
  \includegraphics[width=10cm]{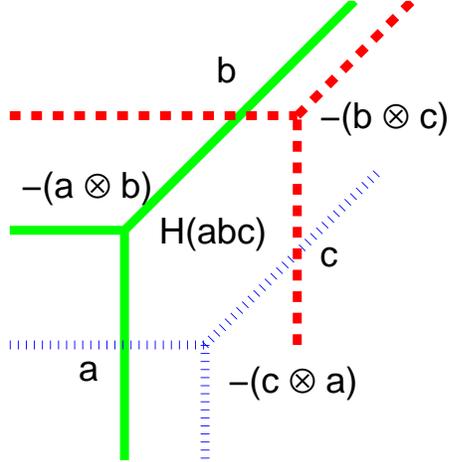}\\
  \caption{Transversal tropical triangle.}
  \label{fig:trian}\end{figure}

\begin{dfn}\label{dfn:transversal_triangle} Three points  $a,b,c$   define a \emph{transversal tropical triangle}  $abc$ if the vertices $a,b,c$ are pairwise transversal and so are the tropical sides.
\end{dfn}

\begin{dfn}\label{dfn:good_triangle}
Three  points $a,b,c$   define a \emph{good tropical triangle}  $abc$ if, by stable join, they give rise to three
tropical lines  $ab, bc, ca$ which,
 stably intersected by pairs, yield the original points $a,b,c$, i.e., $ca\cap_{est}ab= a, ab\cap_{est}bc= b, bc\cap_{est}ca= c$.
\end{dfn}

\begin{dfn}\label{dfn:proper_triangle}
Three points  $a,b,c$   define a \emph{proper tropical triangle}  $abc$ if $a,b,c,-(a\otimes b),-(b\otimes c)$ and $-(c\otimes a)$ are six
different points. A tropical triangle  which is not proper will be called
\emph{improper}.
\end{dfn}

\begin{ex}\label{ex:good_non_prop}
The points $[0,0,0], [1,1,0], [0,1,0]$ define a good tropical triangle and so do the points $[0,0,0], [1,1,0], [1,0,0]$. Both triangles are improper.
The triangle in example \ref{ex:good_and_prop} is good and proper.
\end{ex}

By lemmas \ref{lem:transversal} and \ref{lem:long_expr} and corollary \ref{cor:tres_val},
a  tropical triangle $T$ is transversal if and only if $T$ is good and proper.

\begin{thm}\label{thm:triangle_affine}
Three   points $a=\ptoo a$,
 $b=\ptoo b$, $c=\ptoo c$
in $\R^2$  determine a  transversal tropical triangle  if and only if, perhaps after relabeling, the
following inequalities hold:
    $$a_1<b_1<c_1,\qquad a_2<c_2<b_2,$$
    $$b_1-b_2<a_1-a_2<c_1-c_2.$$
    In particular, these inequalities determine an open unbounded polyhedron in $\R^6$, which can be viewed projectively in $\T\P^5$.
\end{thm}

\begin{proof}
The inequalities hold for  $(a_1,a_2,b_1,b_2,c_1,c_2)\in\R^6$ if and only if they hold for $(a_1+\lambda,a_2+\lambda,b_1+\lambda,b_2+\lambda,c_1+\lambda,c_2+\lambda)$, for $\lambda\in\R$. This proves the last statement.

 Let us assume that the six inequalities hold. Easy computations yield  $a\otimes b=[b_2,b_1,a_1+b_2]$, $b\otimes c=[b_2,c_1,c_1+b_2]$ and $c\otimes a=[c_2,c_1,c_1+a_2]$. Then we obtain $(c\otimes a)\otimes(a\otimes b)=[b_2+a_1+c_1, c_1+a_2+b_2,b_2+c_1]=a$ and, similarly, $(a\otimes b)\otimes(b\otimes c)=b$ and $(b\otimes c)\otimes(c\otimes a)=c$, proving that the tropical triangle $abc$ is transversal.

Conversely, suppose that the points $a,c,b$ define a transversal triangle.  Then no two of the given points lie on a classical line of slope $\infty,0,1$.
In particular, the numbers $a_1-a_2$, $b_1-b_2$, $c_1-c_2$ are pairwise different. We may  assume that $$a_1< \min\{b_1, c_1\}$$ $$b_1-b_2< c_1-c_2.$$ Then we compute the point $(c\otimes a)\otimes(a\otimes b)$ and see that one of its coordinates is given by the value of a singular tropical determinant, unless $$b_1-b_2< a_1-a_2< c_1-c_2.$$
Now classical and tropical geometry tell us that the coordinates of the vertex $-(a\otimes b)$ are $[a_1,a_1-b_1+b_2,0]$, on the one hand, and $[-\max\{a_2,b_2\}, -b_1, -a_1-b_2]$ on the other. Equating these projective tropical  coordinates yields $$\max\{a_2,b_2\}=b_2.$$ A similar computation for  $-(c\otimes a)$ yields $$\max\{a_2,c_2\}=c_2.$$
Therefore, we have 
$$a_2< \min\{b_2, c_2\}.$$
We proceed to determine the values of $\min\{b_1, c_1\}$ and $\min\{b_2, c_2\}$. There are only three possible cases, because the condition $b_1-b_2< c_1-c_2$ eliminates the possibility $b_1>c_1$ and $b_2<c_2$. Now if $b_1>c_1$ and $b_2>c_2$, then
$b\otimes c=[b_2,b_1,c_1+b_2]$ and $a\otimes b=[b_2,b_1,a_1+b_2]$, so that
$(b\otimes c)\otimes(a\otimes b)=-(b\otimes c)\neq b$, contradicting transversality. And if  $b_1<c_1$ and $b_2<c_2$, then
$b\otimes c=[c_2,c_1,c_1+b_2]$ and  $c\otimes a=[c_2,c_1,c_1+a_2]$ so that $(c\otimes a)\otimes(b\otimes c)=-(b\otimes c)\neq c$, contradicting transversality. Therefore,  $b_1<c_1$ and $b_2>c_2$ and all six inequalities have been proved.
\end{proof}

Let us make a picture of a tropical triangle $abc$ in $\R^2$, say in $Z=0$. The six inequalities shown in the previous theorem must be satisfied. We have to represent the  vertices  $a,b,c$  and the tropical sides $ab,bc,ca$, whose vertices are the points
$$-(a\otimes b)=\pto{a_1}{a_1+b_2-b_1},\quad -(c\otimes a)=\pto{a_2+c_1-c_2}{a_2},\quad -(b\otimes c)=\pto{c_1}{b_2}.$$
Then we obtain a \emph{classical convex hexagon}  having vertices  (in clockwise order)
$$a, -(a\otimes b), b, -(b\otimes c),c, -(c\otimes a)$$ and slopes $\infty,1,0,\infty,1,0$, see figure \ref{fig:trian}. It will be denoted $\HH(abc)$. The six inequalities shown in theorem \ref{thm:triangle_affine}  provide the lengths of the sides of $\HH(abc)$.

\m Up to translation, scaling and exchange of variables $X,Y,Z$, a transversal tropical triangle $T$ is determined  by a classical convex hexagon $\HH\subset \R^2$ of slopes $\infty,0,1$. Now,  $\HH$  is determined by the lattice lengths of its sides, which are real positive numbers $l_1, l_2, \ldots, l_6$ such that
$$l_{j-2}+l_{j-1}=l_{j+1}+l_{j+2}, \qquad j=1,2\ (\text{or\ } j=1,2,\ldots,6)$$
where subscripts are taken modulo 6. Therefore, the set
$$P:=\{[l_1, l_2, \ldots, l_6]\in\T\P^5: l_j>0 \text{\ and\ } l_{j-2}\odot l_{j-1}=l_{j+1}\odot l_{j+2},\ j=1,2\}$$ is a \emph{parameter space for transversal tropical triangles}. The dimension of $P$ is three and  any positive numbers $l_1, l_2, l_3, l_5\in\R$  such that
$$l_5<\min\{l_1+l_2, l_2+l_3\}$$ determine a unique point in $P$.



\m
 Let us now translate the six inequalities found in
theorem \ref{thm:triangle_affine} to the
projective setting.

\begin{thm}\label{thm:triangle_projective} Let
$a=\ptoop {a'}$, $b=\ptoop {b'}$, $c=\ptoop {c'}\in\T\P^2$
be three different interior points. Then $a,b,c$ determine
a  transversal tropical triangle  if and only if, perhaps after
relabeling, the following six inequalities hold:
$$b'_1-b'_2<a'_1-a'_2<c'_1-c'_2,$$
    $$a'_2-a'_3<c'_2-c'_3<b'_2-b'_3,$$
    $$c'_3-c'_1<b'_3-b'_1<a'_3-a'_1.$$
\end{thm}
\begin{proof} Without loss of generality, we may work in
$Z=0$. Then $a=\ptop {a'_1-a'_3}{a'_2-a'_3}0$, $b=\ptop
{b'_1-b'_3}{b'_2-b'_3}0$, $c=\ptop
{c'_1-c'_3}{c'_2-c'_3}0$. Now,
 notice that the points $\pto {a'_1-a'_3}{a'_2-a'_3}$,  $\pto
{b'_1-b'_3}{b'_2-b'_3}$, $\pto {c'_1-c'_3}{c'_2-c'_3}$
belong to $\R^2$ and satisfy the six inequalities shown in  theorem
\ref{thm:triangle_affine}, where
$a_j=a'_j-a'_3$, $b_j=b'_j-b'_3$, $c_j=c'_j-c'_3$,
$j=1,2$,
 and so we are done.
\end{proof}

Notice the high cyclic symmetry shown in the six inequalities in theorem \ref{thm:triangle_projective}.
Notice also that the letters $a,b,c$ are arranged as a latin square. 

\m
Now we proceed to relate tropical transversal triangles with the notions of tropical span, tropical independence and strong permanent.
The following  definitions are standard.


\begin{dfn}\label{dfn:gener_indep}
In $\T^n$ or $\T\P^n$, let $\suces u1s$ be different
interior points. A point  $u$ is \emph{tropically spanned
 by $\suces u1s$} if it can be written as
    $$u=\lambda_1\odot u_1\oplus\cdots\oplus  \lambda_s\odot
     u_s=\max\{\lambda_1+u_1 ,\ldots, \lambda_s+u_s\},$$
     for some $\suces {\lambda}1s\in\T$, and not all $\lambda_j$
     equal to $-\infty$.

     The points $\suces u1s$ are \emph{tropically independent} if
there does not exist $j\in \{1,\ldots,s\}$ such that $u_j$
is tropically spanned by
$u_1,\ldots,u_{j-1},u_{j+1},\ldots,u_s$.
\end{dfn}
Notice that all points spanned by interior
points are interior.

\m  Tropical span is closer to classical  convexity than to linear, affine or projective span. This is shown in the following theorem, which can be traced back to \cite{Cuninghame_B}, \cite{Develin}, \cite{Joswig}, \cite{Wagneur_F}, \cite{Wagneur_M}.
The classical segment defined by points $a,b$  will be denoted
$\conv(a,b)$.

\begin{thm}\label{thm:hexagono}
Let three   points $a,b,c$ in $\T\P^2$
determine a transversal tropical triangle. Then the following assertions hold.
\begin{enumerate}
\item The points tropically
spanned by  $a,c, b$  are exactly those of the solid hexagon
$\HH(abc)$.
\item The points $a,b,c$  are tropically independent.
\item The matrix $3\times 3$ given by  coordinates of  $a,b,c$ is tropically regular, (i.e., this matrix has a strong permanent).
\end{enumerate}
\end{thm}
\begin{proof} We may work in $Z=0$ and therefore, we may assume, without loss of generality,
that  the coordinates of $a,b,c$ satisfy the six inequalities shown in theorem
\ref{thm:triangle_affine}.

First, we show that the vertices $b,c$ tropically
span the union of the classical  segments  $\conv(b,-(b\otimes c))$ and
$\conv(-(b\otimes c),c)$. We have $b=\ptop {b_1}{b_2}0$, $c=\ptop
{c_1}{c_2}0$ with    $b_1<c_1$ and $c_2<b_2$. Then
$$-(b\otimes c)=[c_1,b_2,0]=b\oplus c.$$
A point $u$ tropically spanned by $b,c$ is $u=\lambda \odot
b\oplus \mu\odot c$. Working in the projective plane, we may assume $\lambda=0$  and $\mu\in\T$. Thus
$$u=[\max\{b_1, \mu+c_1\},
\max\{ b_2, \mu+c_2\} ,\max\{0, \mu \}].$$

Now, if $0\ge \mu$, then $b_2> \mu+c_2$.
\begin{itemize}
\item If $b_1\ge \mu+c_1$, then  $u=b$.

\item If $b_1\le \mu+c_1$, then  $u=\ptop
{\mu+c_1}{b_2}{0}$ and the point
$\pto{\mu+c_1}{b_2}\in\R^2$ runs through the segment
$\conv(b,-(b\otimes c))$.
\end{itemize}

Now, if $0\le \mu$, then $b_1< \mu+c_1$.
\begin{itemize}
\item If $b_2\ge \mu+c_2$, then  $u=\ptop
{\mu+c_1}{b_2}{\mu}=\ptop {c_1}{b_2-\mu}0$
and the point $\pto{c_1}{b_2-\mu}\in\R^2$ runs through the
segment $\conv(-(b\otimes c),c)$.

\item If $b_2\le \mu+c_2$, then  $u=c$.
\end{itemize}
In particular, we have proved that the point   $u=b\oplus \mu\odot c\in \{b,c\}$, whenever  $\mu$ does not belong to the closed interval $[b_1-c_1,b_2-c_2]$.
As a by--product, we have proved that $a$ is not tropically spanned by $b,c$.

In a similar manner, the vertices $a,c$
tropically span the set   $\conv(a,-(c\otimes a))\cup\conv(-(c\otimes a),c)$ and  $a,b$ tropically span
the set $\conv(a,-(a\otimes b))\cup\conv(-(a\otimes b),b)$.
Therefore, all points on the  sides of the classical hexagon $\HH(abc)$  are tropically spanned by  $a,b,c$.
 And, in particular, the tropical independence of $a,b,c$ follows.

Now  we consider a point $t=\ptoo t\in\R^2$ in the interior of the hexagon
$\HH(abc)$. Then we take the (unique)
point $u=\pto {u_1}{t_2}\in\R^2$ on the border of $\HH(abc)$ with $u_1<t_1$ and
the (unique) point $q=\pto {t_1}{q_2}$ on the border with
$q_2<t_2$. Then $t=u\oplus q$ and, since $u$ is tropically
spanned by $a,b$ and $q$ is tropically spanned by
$a,c$, then $t$ is tropically
spanned by $a,b,c$.

To finish up, let $$A=\left(\begin{array}{ccc} a_1&a_2&0\\b_1&b_2&0\\c_1&c_2&0\\
\end{array}\right)$$ and check that $|A|_{trop}=c_1+b_2$,
using the six inequalities. Moreover,  the maximum is attained only once, showing that $A$ is tropically regular.
\end{proof}

The converse to the second and third statements in theorem \ref{thm:hexagono} do not hold.

\begin{ex} The points $a=[0,0,0]$, $b=[3,9,0]$ and $c=[2,1,0]$ are tropically independent but they only satisfy five of the six inequalities in theorem \ref{thm:triangle_affine}. Moreover, the coordinate matrix $A$ is tropically regular and $|A|_{trop}=11$.
\end{ex}


%
%

\m \label{comment:tree_of_tropical_triangles}
Now let $T$ be a
 transversal (i.e., good and proper) tropical triangle  $abc$  (affine or projective).
Good improper  tropical triangles  arise from $T$, by letting two
or more adjacent vertices $a, -(a\otimes b), b, -(b\otimes c),c, -(c\otimes a)$
of $\HH(abc)$ collapse, but keeping $a,b,c$ pairwise different. This means that  the classical  hexagon
$\HH(abc)$ collapses to a $n$--polygon, with $3\le n<6$
sides (of slopes $\infty,0,1$). Equivalently, at most three inequalities in theorems
\ref{thm:triangle_affine} or
\ref{thm:triangle_projective} become
equalities. The reader can easily sketch (say  in $Z=0$),
the 14 existing combinatorial types of improper good tropical triangles  thus
obtained.  He/she can also arrange them into a graph. The
vertices in this graph are  improper good tropical triangles  and two such
triangles  $T'$ and $T''$ in this graph are joined by an edge if
$\HH(T'')$ is obtained from  $\HH(T')$
 by collapsing two consecutive vertices.
 The leaves in this graph correspond either to  span tropical triangles $T$ such that either $\HH(T)$ is a classical pentagon (there are six such leaves)
 or  $\HH(T)$ is a classical triangle (there are two such leaves: $T=abc$ such that
 $a=-(a\otimes b), b=-(b\otimes c), c=-(c\otimes a)$ or  $b=-(a\otimes b), b=-(b\otimes c), a=-(c\otimes a)$).
 The latter are, by the way, like the two triangles shown in example \ref{ex:good_non_prop}.


\end{document}